\documentclass[12pt,leqno]{amsart}
\usepackage{amssymb}
\usepackage{amsmath,amssymb,color}
\numberwithin{equation}{section}

\newcommand{\ep}{\varepsilon}
\newcommand{\la}{\lambda}
\newcommand{\va}{\varphi}
\newcommand{\ppp}{\partial}

\newcommand{\DDD}{\mathcal{D}}
\newcommand{\www}{\widetilde}

\newcommand{\pppa}{\partial_t^{\alpha}}

\newcommand{\ddda}{d_t^{\alpha}}

\newcommand{\NUNU}{\partial_{\nu_A}}

\newcommand{\R}{\mathbb{R}}
\newcommand{\Q}{\mathbb{Q}}

\newcommand{\C}{\mathbb{C}} 
\newcommand{\N}{\mathbb{N}}

\newcommand{\ooo}{\overline}
\newcommand{\OOO}{\Omega}

\newcommand{\sumn}{\sum_{n=1}^{\infty}}
\newcommand{\sumd}{\sum_{i,j=1}^d}

\newcommand{\CC}{{_{0}C^1[0,T]}}

\newcommand{\HH}{H_{\alpha}}
\newcommand{\MLEE}{(t-s)^{\alpha-1}E_{\alpha,\alpha}(-\lambda_n(t-s)^{\alpha})}

\allowdisplaybreaks

\title
[]
{
Uniqueness for inverse source problems for fractional diffusion-wave equations
by data during not acting time
}

\author{
$^{1,2,3}$ M.~Yamamoto }

\thanks{
$^1$ Graduate School of Mathematical Sciences, The University
of Tokyo, Komaba, Meguro, Tokyo 153-8914, Japan \\
$^2$ Honorary Member of Academy of Romanian Scientists, 
Ilfov, nr. 3, Bucuresti, Romania \\
$^3$ Correspondence member of Accademia Peloritana dei Pericolanti,\\
Palazzo Universit\`a, Piazza S. Pugliatti 1 98122 Messina Italy \\
e-mail: {\tt myama@ms.u-tokyo.ac.jp}
}

\date{}
\begin{document}
\maketitle

\baselineskip 18pt

\begin{abstract}
We consider fractional diffusion-wave equations with source term which 
is represented in a form of a product of a temporal function and 
a spatial function.
We prove the uniqueness for 
inveres source problem of determining spatially varying factor 
by decay of data as the time tends to $\infty$, provided that 
the source does not work during the observations.
Our main result asserts the uniqueness if data decay more rapidly than 
$\left(\frac{1}{t^p}\right)$ with any $p\in \N$ as $t\to\infty$.
Date taken not from the initial time are realistic but the uniqueness was not
known in general.
The proof is based on the analyticity and the asymptotic behavior of a function
generated by the solution.
\end{abstract}


\section{Introduction}

Let $\OOO \subset \R^d$ be a bounded domain with smooth boundary
$\ppp\OOO$, $\nu(x) = (\nu_1(x), ..., \nu_d(x))$ be the unit outward normal 
vector to $\ppp\OOO$ at $x$, and let
$$
0< \alpha < 2, \quad \alpha \ne 1.
$$
Throughout this article, we set  
$$
(-Av)(x) = \sumd \ppp_i(a_{ij}(x)\ppp_jv(x)) + c(x)v(x), \quad 
x\in \OOO,                            \eqno{(1.1)}
$$
where $a_{ij} = a_{ji} \in C^1(\ooo{\OOO})$, $1\le i,j\le d$ and 
$c\in C(\ooo{\OOO})$ and $c\le 0$ in $\OOO$.
Moreover we assume that there exists a constant $\kappa>0$ such that 
$$
\sumd a_{ij}(x)\zeta_i\zeta_j \ge \kappa\sum_{i=1}^d \zeta_i^2 
\quad \mbox{for all $x \in \ooo{\OOO}$ and $\zeta_1, ..., \zeta_d \in \R$}.   
$$
We set 
$$
\NUNU v(x) := \sumd a_{ij}(x)(\ppp_jv)\nu_i.
$$

We consider an initial boundary value problem for a fractional
diffusion-wave equation whose source term is described by 
$\mu(t)f(x)$, where $f$ is a spatial distribution of the source 
and $\mu$ is a temporally changing factor.  We can formally describe
an initial boundary value problem:
$$
\left\{ \begin{array}{rl}
& \ddda u(x,t) = -Au(x,t) + \mu(t)f(x), \quad x\in \OOO, \, 0<t<T, \\
& u\vert_{\ppp\OOO\times (0,T)} = 0, \\
& u(x,0) = 0, \quad x \in \OOO \quad \mbox{if $0<\alpha < 1$}, \\
& u(x,0) = \ppp_tu(x,0) = 0, \quad x \in \OOO 
\quad \mbox{if $1<\alpha < 2$}.
\end{array}\right.
                                       \eqno{(1.2)}
$$
Here we can define the pointwise classical Caputo derivative:
$$
\ddda v(t) = 
\left\{ \begin{array}{rl}
& \frac{1}{\Gamma(1-\alpha)}\int^t_0
(t-s)^{-\alpha}\frac{d}{ds}v(s) ds, \quad 0<\alpha<1
\quad \mbox{for $v \in W^{1,1}(0,T)$}, \\
& \frac{1}{\Gamma(2-\alpha)}\int^t_0
(t-s)^{1-\alpha}\frac{d^2}{ds^2}v(s) ds, \quad 1<\alpha<2
\quad \mbox{for $v \in W^{2,1}(0,T)$}.
\end{array}\right.
$$

We postpone a rigorous formulation of the initial boundary value problem,
and we first describe our motivation for the inverse source problem
in this article.

The first equation in (1.2) describes anomalous diffusion in 
heterogeneous media such as soil.  In the heterogenous media, it is often
observed that the classical diffusion equation corresponding to 
$\alpha=1$ cannot simulate diffusion phenomenon well and real density 
profiles are deviated from the ones given by the equation with $\alpha=1$.
Therefore the time-fractional diffusion-wave equations in (1.2) with 
$\alpha  \in (0,1) \cup (1, 2)$ have been proposed and studied
(e.g., Metzler and Klafter \cite{MK}).  
The right-hand side $\mu(t)f(x)$ of the fractional equation in (1.2) models 
a source term causing diffusion of some substances, for example, 
harmful contaminants such as cesium-137 after destruction of plants or 
storages.  The source of the diffusion can be often modelled by a 
product of a temporal function $\mu(t)$ and a spatial density 
profile $f(x)$.  In a case of radioactive isotope, we can choose
$\mu(t) = e^{-\mu_0t}$ with some decay constant $\mu_0>0$.

Usually we cannot a priori know the source term, for example, when
some dangerous substances start to diffuse 
after some incident in an inaccesible place.
Then a serious direct problem is concerned with accurate 
prediction of contaminants, and for more realistic
simulations, we have to estimate or identify $f(x)$ and/or
$\mu(t)$ in the source by observation data. 

Thus, in this article, we are concerned with an inverse source problem
of determining $f(x)$, $x\in \OOO$ or $\mu(t)$ in some time interval
by asymptotic behavior of data $u(\cdot,t)\vert_{\omega}$ and
$\ppp_{\nu_A}u(\cdot,t)\vert_{\gamma}$ as $t\to\infty$ of  
the solution $u(x,t)$ to (1.2).
Here $\omega \subset \OOO$ is a non-empty subdomain and
$\gamma \subset \ppp\OOO$ is a non-empty subboundary.

Our main interest is the uniqueness in the inverse source problem:
{\it Which behavior of data (1.3) or (1.4) as $t\to \infty$, can 
uniquely determine $f(x)$ for $x\in \OOO$
or $\mu(t)$ over some time interval?}

For the cases where observations start at $t=0$, 
there are many works available on the 
stability as well as the uniqueness.  Here, owing to the numerous 
growing references, we are restricted to a very small number of 
works: Kian, Liu and Yamamoto \cite{KLY}, Kian, Soccorsi, Xue and 
Yamamoto \cite{KSXY}, and a survey Liu, Li and Yamamoto \cite{LLY}.
The readers can consult the references in the above articles.
For the cases of $\alpha=1$ and $=2$, we can refer to 
Yamamoto \cite{Y1}, \cite{Y2}, \cite{Y3}, \cite{Y4}.

By taking practical situations into consideration, it is not realistic
that we can start observations at time $t=0$.  In the case where the diffusion 
was caused by some drastic changes such as explosion, only after 
the incident, one can start to take observation data, and it is more acceptable
to observe data after time passes sufficiently long.  
This is our primary interest for data, and in this article, as is
formulated by (1.7) and (1.8), we adopt 
the asymptotic behavior of solution limited to a subdomain or 
a subboundary at $t\to \infty$.  
For inverse source problems, there are very few theoretical works 
by data measured after positive time.  See Cheng, Lu and Yamamoto 
\cite{CLY} for $\alpha=1,2$.  For $\alpha\in (0,1) \cup (1,2)$, we can 
refer to Theorem 2.2 in Kian, Liu and Yamamoto \cite{KLY},
which proved the uniqueness with
data $u\vert_{\omega \times (t_0,T_0)}$ if $\mu \not\equiv 0$ 
but we have to extra assume that $f=0$ in $\omega$.

To the best knowledge of the author, without such an extra condition on $f$,
there are no published results for the uniqueness of the inverse source 
problem.
\\

For the statement of the main result, we need to formulate an initial boundary 
value problem and the class of solutions where we do not require 
the differentiability for the sake of our choice of natural regularity:
$\mu f \in L^2(0,T;L^2(\OOO))$.
First, we extend the domain of the classical Caputo derivative.

For it, we define the Riemann-Liouville fractional integral operator 
$J^\alpha$ for $\alpha>0$ by 
$$
J^\alpha w(t):= \frac{1}{\Gamma(\alpha)}
\int_0^t (t-s)^{\alpha-1}w(s) ds,
\quad w\in L^2(0,T).
$$
We set 
$$
\CC := \{ w \in C^1[0,T];\, w(0) = 0\}.
$$

Henceforth $\ooo{B}^X$ denotes the closure of a set $B$ in the space
$X$.

Recalling the Sobolev-Slobodecki space and 
$\Vert \cdot \Vert_{H^{\alpha}(0,T)}$ for $0<\alpha<1$
(e.g., Adams \cite{Ad}), we define
$$
\HH(0,T) = \ooo{\CC}^{H^{\alpha}(0,T)}.
$$
It is known (e.g., Gorenflo, Luchko and Yamamoto \cite{GLY},
Kubica, Ryszewska and Yamamoto \cite{KRY}) that 
$J^{\alpha}: L^2(0,T) \longrightarrow H_{\alpha}(0,T)$ 
is isomorphism for $0<\alpha<1$.
Then we define
$$
\pppa := (J^{\alpha})^{-1}, \quad \DDD(\pppa) = \HH(0,T).
$$
Next we define $\pppa$ for $1<\alpha<2$.
Let $\alpha := 1 + \beta$ with $0< \beta \le 1$ and
$$
H_1(0,T) := \{ v \in H^1(0,T);\, v(0) = 0\}
$$
$$
H_{1+\beta}(0,T) := \left\{ v \in H_1(0,T);\, 
\frac{dv}{dt}\in H_{\beta}(0,T)\right\}
$$
with 
$$
\Vert v\Vert_{H_{1+\beta}(0,T)}
:= \Vert v\Vert_{H^1(0,T)} + \left\Vert \frac{dv}{dt}\right\Vert
_{H_{\beta}(0,T)}.
$$
By \cite{GLY} or \cite{KRY}, we can prove that 
$J^{1+\beta}: L^2(0,T) \longrightarrow H_{1+\beta}(0,T)$ is 
isomorphism, and we define
$$
\pppa = \ppp_t^{1+\beta} = (J^{\alpha})^{-1}, \quad  
\mathcal{D}(\pppa) = H_{\alpha}(0,T).
$$
We can directly see that $\pppa v = \ddda v$ for 
$v\in C^2[0,T]$ satisfying $v(0) = \frac{dv}{dt}(0) = 0$.
Throughout this article, based on this extended fractional derivative $\pppa$
with $\alpha \in (0,1) \cup (1,2)$, we consider the initial boundary value 
problem
$$
\left\{ \begin{array}{rl}
& \pppa u(x,t) = -Au(x,t) + \mu(t)f(x), \quad x\in \OOO, \, 0<t<T, \\
& u\vert_{\ppp\OOO\times (0,T)} = 0, \\
& u(x,\cdot) \in \HH(0,T) \quad \mbox{for almost all $x\in \OOO$.}
\end{array}\right.
                                       \eqno{(1.3)}
$$
If $\frac{1}{2} < \alpha$ and $u(x,\cdot) \in H_{\alpha}(0,T)$, 
then Sobolev embedding $u(x,\cdot) \in C[0,T]$, and $u(x,\cdot)$ can be 
approximated by functions in $\CC$ and so $u(x,\cdot) \in \HH(0,T)$
is interpreted as $u(x,t) \longrightarrow 0$ as $t \to 0$ for almost all
$x\in \OOO$.  In the case $\frac{3}{2} < \alpha$, we can similarly
interpret the zero initial condition for $u(\cdot,0)$ and
$\ppp_tu(\cdot,0)$ within the framework $\HH(0,T)$. 
Thus the third condition in (1.5) means the generalized zero initial 
condition in a unified manner for each $\alpha \in (0,1) \cup (1,2)$.
Then we can prove that for $f\in H^2(\OOO)\cap H^1_0(\OOO)$ 
and $\mu\in L^2(0,T)$, there 
exists a unique solution 
$$
u \in L^2(0,T;H^2(\OOO)\cap H^1_0(\OOO))\cap \HH(0,T;L^2(\OOO))
                                                  \eqno{(1.4)}
$$
to (1.3).

We can relax the regularity $H^2(\OOO) \cap H^1_0(\OOO)$ of $f$ for (1.4), 
but we do not pursue here.
The proof of (1.4) is found in \cite{KRY} for $0<\alpha<1$ and
in Theorem 2.2 (ii) in \cite{SY} for $1<\alpha<2$, for example.  Indeed 
for $0<\alpha<1$, the work \cite{KRY} establishes the unique 
existence of $u$ satisfying (1.4) under regulairty 
$f \in L^2(\OOO)$.  Furthermore we can weaken the condition 
$\mu \in L^2(0,T)$, but we do not discuss details.

Throughout this article, we always consider the solution class (1.4) for 
(1.3).

Let $t_0>0$ be arbitrarily fixed.
We assume 
$$
\mu \in L_{loc}^2(0,\infty), 
\quad \mu(t) = 0 \quad \mbox{for $t>t_0$}.             
                                                      \eqno{(1.5)}
$$

Let $\Vert \cdot\Vert$ and $(\cdot, \cdot)$ denote the norm and the 
scalar product in $L^2(\OOO)$.
For an elliptic operator $A$ defined by (1.1), we define the domain 
$\DDD(A)$ by $\DDD(A) = H^2(\OOO) \cap H^1_0(\OOO)$.  Then the spectrum 
$\sigma(A)$ of 
$A$ is composed of positive eigenvalues with finite multiplicities.
The set $\sigma(A) = \{ \la_n\}_{n\in \N}$ is numbered as 
$$
0 < \la_1 < \la_2 < \cdots \longrightarrow \infty.
$$
Let $\{ \va_{nj}\}_{1\le j\le d_n}$ be an orthonormal basis of Ker $(\la_n-A)$.
We further set
$$
P_na = \sum_{j=1}^{d_n} (a, \va_{nj})\va_{nj}, \quad n\in \N, \, 
a\in L^2(\OOO).
$$
Then $a = \sum_{n=1}^{\infty} P_na$ in $L^2(\OOO)$ for all $a\in L^2(\OOO)$.
Moreover we can define a fractional power $A^{\sigma}$ with $\sigma>0$
(e.g., Pazy \cite{Pa}) and we note 
$\Vert v\Vert_{H^{2\sigma}(\OOO)} \le C\Vert A^{\sigma}v\Vert$ for 
each $v\in \DDD(A^{\sigma})$.

Let two numbers $\sigma_1$ and $\sigma_2$ satisfy 
$$
\sigma_1 > \max \left\{ 1, \, \frac{d}{2}-2\right\}, \quad
\sigma_2 > \max \left\{ 1, \, \frac{d}{2}-\frac{1}{4}\right\}.
                                                            \eqno{(1.6)}
$$

Now we are ready to state the main result.
\\
{\bf Theorem 1.}
\\
{\it
Let $\alpha \in (0,1) \cup (1,2)$ and let (1.5) hold.
\\
(i) Let $f\in \DDD(A^{\sigma_1})$.  
We assume that for any $p\in \N$, there exists a constant
$C(p) > 0$ such that 
$$
\Vert u(\cdot,t)\Vert_{L^2(\omega)} \le \frac{C(p)}{t^p}
\quad \mbox{as $t\to \infty$}.                     \eqno{(1.7)}
$$
Then $\mu=0$ in $(0,t_0)$ or $f=0$ in $\OOO$.
\\
(ii) Let $f\in \DDD(A^{\sigma_2})$.  
We assume that for any $p\in \N$, there exists a constant
$C(p) > 0$ such that 
$$
\Vert \ppp_{\nu_A}u(\cdot,t)\Vert_{L^2(\gamma)} \le \frac{C(p)}{t^p}
\quad \mbox{as $t\to \infty$}.                     \eqno{(1.8)}
$$
Then $\mu=0$ in $(0,t_0)$ or $f=0$ in $\OOO$.
}
\\

In general, as is seen by the arguments in Section 2 (especially
by e.g., (2.7)), we can find a constant $\beta > 0$ such that 

$$
\Vert u(\cdot,t)\Vert_{L^2(\omega)} = O\left( \frac{1}{t^{\beta}}\right),
\quad
\Vert \ppp_{\nu_A}u(\cdot,t)\Vert_{L^2(\omega)} 
= O\left( \frac{1}{t^{\beta}}\right)
$$
as $t\to \infty$.  Thus the assumptions (1.7) and (1.8) especially require
more rapid decay, in order to conclude that $u=0$ in $\OOO\times (0,\infty)$.
\\

By means of the analyticity of $u(\cdot,t)$ for $t>t_0$, from Theorem 1
we can derive
\\
{\bf Corollary 1.}
\\
{\it
We arbitrarily choose $T_0 > t_0$.
Let $\alpha \in (0,1) \cup (1,2)$ and let (1.5) hold.
\\
(i) Let $f\in \DDD(A^{\sigma_1})$.  
Then $u=0$ in $\omega \times (t_0,T_0)$ implies 
$\mu=0$ in $(0,t_0)$ or $f=0$ in $\OOO$.
\\
(ii) Let $f\in \DDD(A^{\sigma_2})$.  
Then $\ppp_{\nu_A}u = 0$ on $\gamma \times (t_0,T_0)$ implies 
$\mu=0$ in $(0,t_0)$ or $f=0$ in $\OOO$.
}
\\

This article is composed of five sections.
In Section 2, we construct $t$-functions related to the solution $u$ 
and prove Proposition 1 which is a key for the proof of Theorem 1.
In Section 3, we complete the proof of Theorem 1 and in Section 4 we provide
the proof of Corollary 1.  Finally Section 5 gives
concluding remarks.
\section{Key proposition}

For $\alpha, \beta > 0$, we define the Mittag-Leffler functions by 
$$
E_{\alpha,\beta}(z) = \sum_{k=0}^{\infty} \frac{z^k}{\Gamma(\alpha k+ \beta)},
$$  
which is an entire function in $z\in \C$.

We prove
\\
{\bf Lemma 1.}
\\
{\it
Let $T>0$ be arbitrary.
For $\mu \in L^2(0,T)$, the unique solution
$u$ in the class (1.6) to (1.5) is given by 
$$
u(x,t) = \sumn \left( \int^t_0 \MLEE \mu(s) ds\right) (P_nf)(x)
                                                     \eqno{(2.1)}
$$
in $\HH(0,T;L^2(\OOO)) \cap L^2(0,T;H^2(\OOO) \cap H^1_0(\OOO))$.
}
\\

The proof relies on the classical Fourier method by the separation of
variables, or is found for example in Sakamoto and Yamamoto \cite{SY}.

We choose $T_1>0$ such that 
$$
t_0 < T_1
$$
and set 
$$
\psi_n(t) := \int^{t_0}_0 \MLEE \mu(s) ds, \quad t\ge T_1, \quad 
n\in \N.                         \eqno{(2.2)}
$$

Henceforth $C>0$, $C_k>0$ denote generic constants which are independent of
$n\in \N$, $f$ and choices of $t$ in intervals under consideration.  
When we emphasize the dependence on other quantities $K$, we describe
$C(K)$, $C_1(K)$, etc., but we omit the dependence if there is no fear of
confusion.

Then
\\
{\bf Lemma 2.}
\\
{\it
We see $\psi_n \in C[T_1,\, \infty)$ and there exists a constant
$C>0$ such that 
$$
\vert \psi_n(t)\vert \le \frac{C}{\la_n}\Vert \mu\Vert_{L^1(0,t_0)}
                                                        \eqno{(2.3)}
$$
for all $n\in \N$ and $t \ge T_1$.
}
\\
{\bf Proof.}
\\
Since $t \ge T_1 > t_0$, we see that $0<T_1-t_0 < t-s$ for all 
$0<s\le t_0$, the function $\MLEE$ is continous in $t \ge T_1$ and 
$0\le s \le t_0$, so that $\psi_n\in C[T_1,\, \infty)$.
Morever, by 
$$
\vert E_{\alpha,\alpha}(-\eta)\vert \le \frac{C}{1+\eta}, \quad 
\eta > 0, \quad 0<\alpha < 2
$$
(e.g., Theorem 1.6 (p.35) in Podlubny \cite{Po}), using 
$\vert t-s\vert \ge T_1-t_0>0$, we see 
$$
\vert \psi_n(t)\vert \le \int^{t_0}_0 (t-s)^{\alpha-1}
\frac{C\vert \mu(s)\vert}{1+ \la_n\vert t-s\vert^{\alpha}} ds.
$$
Here we note that $\la_n>0$ for all $n\in \N$.
Thus the proof of Lemma 2 is complete.
$\blacksquare$
\\

Now, we show the key proposition.
\\
{\bf Proposition 1.}
\\
{\it 
Let $\mu \not\equiv 0$ in $(0,t_0)$.
We assume that $a_n\in \R$, $n\in \N$ satisfy 
$$
\sumn \left\vert \frac{a_n}{\la_n} \right\vert < \infty.    \eqno{(2.4)}
$$
If for each $p\in \N$, there exist constants $T_2(p)>0$ and $C(p) > 0$
such that 
$$
\left\vert \sumn a_n\psi_n(t)\right\vert \le \frac{C(p)}{t^p},
\quad t > T_2(p),               \eqno{(2.5)}
$$
then $a_n=0$ for all $n\in \N$.
}
\\

By Lemma 3 and (2.4), we notice that the series $\sumn a_n\psi_n$ is uniformly 
convergent in any compact set in $(T_1, \,\infty)$.
\\
{\bf Proof of Proposition 1.}
\\
We first prove
\\
{\bf Lemma 3.}\\
{\it 
Let (2.4) hold.  We assume that there exists $\{\ell_k\}_{k\in \N} 
\subset \N$ satisfying $\lim_{k\to \infty} \ell_k=\infty$.
If 
$$
\sumn \frac{a_n}{\la_n^{\ell_k}} = 0 \quad \mbox{for all $k\in \N$},
$$
then $a_n=0$ for all $n\in \N$.
}
\\
{\bf Proof of Lemma 3.}\\
We have
$$
\frac{a_1}{\la_1^{\ell_k}} + \sum_{n=2}^{\infty} 
\frac{a_n}{\la_n^{\ell_k}} = 0,\quad \mbox{that is},
\quad 
a_1 + \sum_{n=2}^{\infty} a_n \left( \frac{\la_1}{\la_n}
\right)^{\ell_k} = 0.
$$
Hence 
$$
\vert a_1\vert = \left\vert -\sum_{n=2}^{\infty} a_n 
\left( \frac{\la_1}{\la_n}\right)^{\ell_k}\right\vert
= \left\vert \sum_{n=2}^{\infty} 
\left( a_n \frac{\la_1}{\la_n}\right)
\left( \frac{\la_1}{\la_n}\right)^{\ell_k-1} \right\vert
\le \vert \la_1\vert  
\left( \sum_{n=2}^{\infty} \left\vert \frac{a_n}{\la_n}\right\vert \right)
 \left( \frac{\la_1}{\la_2}\right)^{\ell_k-1}.
$$
By $0 < \la_1 < \la_2 < ....$, we see that 
$\left\vert \frac{\la_1}{\la_2}\right\vert < 1$.
Letting $k \to \infty$, we see that $\ell_k \to \infty$, and so
$a_1 = 0$.  Therefore, 
$$
\sum_{n=2}^{\infty} \frac{a_n}{\la_n^{\ell_k}} = 0, \quad k\in \N.
$$
Repeating the above argument, we have $a_2=a_3= \cdots = 0$.
Thus the proof of Lemma 3 is complete.
$\blacksquare$

We note the asymptotic expansion of the Mittag-Leffler function:
$$
\eta^{\alpha-1}E_{\alpha,\alpha}(-\la\eta^{\alpha})
= \sum_{\ell=2}^N \frac{(-1)^{\ell+1}}{\Gamma(\alpha-\alpha\ell)}
\frac{1}{\la^{\ell}\eta^{\alpha\ell-\alpha+1}}
+ O\left( \frac{1}{\la^{N+1}\eta^{\alpha(N+1)-\alpha+1}}\right)
$$
as $\eta \to \infty$ for each $N\in \N \setminus \{1\}$
(e.g., Theorem 1.4 (pp.33-34) in \cite{Po}).
We arrange a set $\{ \ell\in \N\setminus \{1\};\, 
\alpha\ell \not\in \N\}$ in an increasing
manner:
$$
\ell_1 < \ell_2 < \cdots. 
$$
Here $\{ \ell\in \N\setminus \{1\};\, \alpha\ell \not\in \N\}$ 
is an infinite set, and so 
$\lim_{k\to\infty} \ell_k = \infty$.
Indeed, if it is not a finite set, then we see by 
$\N \setminus \{1\}
= \{ \ell\in \N\setminus \{1\};\, \alpha\ell \not\in \N\} \cap 
\{ \ell\in \N\setminus \{1\};\, \alpha\ell \in \N\}$ that there exists 
some $N_0 \in \N$ such that 
$$
\{ \ell\in \N\setminus \{1\};\, \alpha\ell \in \N\} 
\supset \{N_0, N_0+1, N_0+2, \cdots\},
$$
so that $\alpha(N_0+1),\, \alpha N_0\in \N$.  Therefore, 
$\alpha = \alpha(N_0+1) - \alpha N_0 \in\N$, which is impossible by 
$\alpha \in (0,1) \cup (1,2)$.  $\blacksquare$
\\

In particular, for $\alpha \not\in \Q$, we note that $\ell_k = k$ for
$k=2,3,4,...$.

Hence, for each $K\in \N$, we have
$$
\eta^{\alpha-1}E_{\alpha,\alpha}(-\la\eta^{\alpha})
= \sum_{k=1}^K \frac{(-1)^{\ell_k+1}}{\Gamma(\alpha-\alpha\ell_k)}
\frac{1}{\la^{\ell_k}\eta^{\alpha\ell_k-\alpha+1}}
+ O\left( \frac{1}{\la^{\ell_{K+1}}\eta^{\alpha\ell_{K+1}-\alpha+1}} 
\right)
$$
as $\la\eta^{\alpha} \to \infty$.
Therefore,
$$
\eta^{\alpha-1}E_{\alpha,\alpha}(-\la_n\eta^{\alpha})
= \sum_{k=1}^K \frac{(-1)^{\ell_k+1}}{\Gamma(\alpha-\alpha\ell_k)}
\frac{1}{\la_n^{\ell_k}\eta^{\alpha\ell_k-\alpha+1}}
+ O\left( \frac{1}{\la_n^{\ell_{K+1}}\eta^{\alpha\ell_{K+1}-\alpha+1}} 
\right)
$$
as $\la_n\eta^{\alpha} \to \infty$ for each $K\in \N$ and $n\in \N$.
Using $\la_n \ge \la_1 > 0$ for $n\in \N$, we see: for each $K\in \N$, 
there exists some constant $\eta_0 = \eta_0(K) > 0$ and 
$C=C(K)>0$ such that if $\la_1\eta^{\alpha} \ge \eta_0(K)$, then 
$$
\left\vert \eta^{\alpha-1}E_{\alpha,\alpha}(-\la_n\eta^{\alpha})
- \sum_{k=1}^K \frac{(-1)^{\ell_k+1}}{\Gamma(\alpha-\alpha\ell_k)}
\frac{1}{\la_n^{\ell_k}\eta^{\alpha\ell_k-\alpha+1}} \right\vert
\le \frac{C(K)}{\la_n^{\ell_{K+1}}\eta^{\alpha\ell_{K+1}-\alpha+1}}
                                          \eqno{(2.6)}
$$
for all $n\in \N$.

Here we note that the constant $C(K)>0$ is independent of $n\in \N$.
\\

We fix $K\in \N$ arbitrarily.
Let $t \ge T_1$.
Since $\la_1(t-s)^{\alpha} \ge \la_1(t-t_0)^{\alpha}$ for 
$0\le s\le t_0$, if $\la_1(t-t_0)^{\alpha} \ge \eta_0$, then 
$\la_1(t-s)^{\alpha} \ge \eta_0$ for $0\le s\le t_0$.
Then, we substitute (2.6) into (2.2), and we obtain
$$
\left\vert 
\psi_n(t) - \sum_{k=1}^K \frac{(-1)^{\ell_k+1}}{\Gamma(\alpha-\alpha\ell_k)}
\frac{1}{\la_n^{\ell_k}}\int^{t_0}_0 
\frac{\mu(s)}{(t-s)^{\alpha\ell_k-\alpha+1}} ds\right\vert
$$
$$
\le \int^{t_0}_0 \frac{C_1}{\la_n^{\ell_{K+1}}}
(t-s)^{\alpha\ell_{K+1}-\alpha+1}\vert \mu(s)\vert ds
$$
for $\la_1(t-t_0)^{\alpha} \ge \eta_0$ and $t>T_1$.

We set
$$
T_3:= \max\left\{ \left( \frac{\eta_0}{\la_1}\right)^{\frac{1}{\alpha}}
+ t_0, \, T_1\right\}.
$$

Multiplying by $a_n$ and taking $\sumn$, by the triangle inequality, we obtain
\begin{align*}
& \left\vert \sumn a_n\psi_n(t) 
- \sumn a_n \sum_{k=1}^K \frac{(-1)^{\ell_k+1}}{\Gamma(\alpha-\alpha\ell_k)}
\frac{1}{\la_n^{\ell_k}}\int^{t_0}_0 
\frac{\mu(s)}{(t-s)^{\alpha\ell_k-\alpha+1}} ds\right\vert\\
\le & C_1\sumn \frac{\vert a_n\vert}{\la_n^{\ell_{K+1}}}
\int^{t_0}_0 \frac{\vert \mu(s)\vert}
{(t-s)^{\alpha\ell_{K+1}-\alpha+1}} ds
\end{align*}
for $t > T_3$.

We write
$$
A_k:= \sumn \frac{a_n}{\la_n^{\ell_k}}, \quad k\in \N, \quad
T^*:= \max\{ T_1, \, T_2(p), \, T_3\}.
$$
Henceforth we omit the dependence of the constants $C>0$, $C_2>0$, etc. on 
$K, \ell_k$, etc.

Then, since $\sumn a_n\psi_n(t) = 0$ for $t>T_1$ and 
$$
\sumn \frac{\vert a_n\vert}{\la_n^{\ell_{K+1}}}
\le C\sumn \left\vert \frac{a_n}{\la_n}\right\vert < \infty,
$$
and 
$$
\int^{t_0}_0 \frac{\vert \mu(s)\vert}{(t-s)^{\alpha\ell_{K+1}-\alpha+1}} ds
\le \frac{\Vert \mu\Vert_{L^1(0,t_0)}}{(t-t_0)^{\alpha\ell_{K+1}-\alpha+1}}
\le \frac{C_2}{(t-t_0)^{\alpha\ell_{K+1}-\alpha+1}},
$$
we obtain
$$
\left\vert \sumn a_n\psi_n(t)
- \sum_{k=1}^K \frac{(-1)^{\ell_k+1}}{\Gamma(\alpha-\alpha\ell_k)}
\left( \int^{t_0}_0 
\frac{\mu(s)}{(t-s)^{\alpha\ell_k-\alpha+1}} ds\right)A_k \right\vert
\le \frac{C_3}
{(t-t_0)^{\alpha\ell_{K+1}-\alpha+1}}
                                                \eqno{(2.7)}
$$
for all $t > T^*$.
Thus, in view of (2.5), we reach
$$
\left\vert  
\sum_{k=1}^K \frac{(-1)^{\ell_k+1}}{\Gamma(\alpha-\alpha\ell_k)}
\left( \int^{t_0}_0 
\frac{\mu(s)}{(t-s)^{\alpha\ell_k-\alpha+1}} ds\right) A_k \right\vert
\le  \frac{C_3}{(t-t_0)^{\alpha\ell_{K+1}-\alpha+1}} 
+ \frac{C(p)}{t^p}
                                                \eqno{(2.8)}
$$
for all $t>T^*$.

On the other hand, we consider the binomial asymptotic expansion.
We fix $0<r<1$ arbitrarily. 
Then, for each $M\in \N$ and $\sigma \in \R$, we know  
$$
(1-\eta)^{-\sigma} = \sum_{m=0}^M 
\left(
\begin{array}{c}
-\sigma \\
m \\
\end{array}
\right) (-\eta)^m + R_{\sigma,M}(\eta), \quad \vert \eta \vert \le r,
                                                                  \eqno{(2.9)}
$$
where we set 
$$
\left(
\begin{array}{c}
-\sigma \\
m \\
\end{array}
\right) := \frac{(-\sigma)(-\sigma-1) \cdots (-\sigma-m+1)}{m!},
$$
and the function $R_{\sigma,M}$ in $\eta$ satisfies 
$$
\vert R_{\sigma,M}(\eta)\vert \le C(M,\sigma)\vert \eta\vert^{M+1}
\quad \mbox{for $\vert \eta\vert \le r$}
$$
with a constant $C(M,\sigma)>0$.  Whenever we fix the ranges of 
$\sigma\in \R$
and $M \in \N$ among finite sets, we can write $C(M,\sigma)$ 
simply by $C_4$.

Henceforth the variables $s$ and $t$ satisfy 
$$
0\le s \le t_0, \quad t>T^*.             \eqno{(2.10)}
$$
Then, by $t>T_1$, we have
$$
\frac{s}{t} \le \frac{t_0}{T_1}.      
$$
Since $t_0 < T_1$, setting $r:= \frac{t_0}{T_1}$, we see that 
$r < 1$, that is,
$$
\frac{s}{t} \le \frac{t_0}{T_1} = r < 1.                \eqno{(2.11)}
$$
In terms of (2.11), applying (2.9), we obtain
$$
\frac{1}{(t-s)^{\sigma}} = \frac{1}{t^{\sigma}}\frac{1}{\left(
1- \frac{s}{t}\right)^{\sigma}}
= t^{-\sigma}\left( 1 - \frac{s}{t}\right)^{-\sigma}
$$
$$
= \sum_{m=0}^M   
\left(
\begin{array}{c}
-\sigma \\
m \\
\end{array}
\right) t^{-\sigma-m}(-s)^m + R_M(t,s), \quad 
0\le s\le t_0,\, t>T^*,                      \eqno{(2.12)}
$$
where 
$$
\vert R_M(t,s)\vert \le C_4t^{-\sigma-M-1}t_0^{M+1}.
$$ 
We further set 
$$
\mu_m:= \int^{t_0}_0 (-s)^m \mu(s) ds, \quad m\in \N \cup \{ 0\}.
                                                                \eqno{(2.13)}
$$
Therefore, in view of (2.12), we see
$$
\int^{t_0}_0 \frac{\mu(s)}{(t-s)^{\sigma}} ds
= \sum_{m=0}^M 
\left(
\begin{array}{c}
-\sigma \\
m \\
\end{array}
\right) \frac{1}{t^{\sigma+m}}\int^{t_0}_0 (-s)^m\mu(s) ds
+ \int^{t_0}_0 \mu(s) R_M(t,s)ds
$$
$$
= \sum_{m=0}^M 
\left(
\begin{array}{c}
-\sigma \\
m \\
\end{array}
\right) \frac{\mu_m}{t^{\sigma+m}}
+ \int^{t_0}_0 \mu(s) R_M(t,s)ds.                         \eqno{(2.14)}    
$$

Consequently 
$$
\left\vert \int^{t_0}_0 \mu(s)R_M(t,s) ds\right\vert 
\le \frac{C_4t_0^{M+1}}{t^{\sigma+M+1}}\Vert \mu\Vert_{L^1(0,t_0)}
\le \frac{C_5}{t^{\sigma+M+1}}, \quad t>T^*.
$$
Hence, in terms of (2.14), by (2.8), we have
$$
\biggl\vert  
\sum_{k=1}^K \frac{(-1)^{\ell_k+1}}{\Gamma(\alpha-\alpha\ell_k)}
A_k \biggl[ \sum_{m=0}^M 
\left(
\begin{array}{c}
-\alpha\ell_k+\alpha-1 \\
m \\
\end{array}
\right)
\frac{\mu_m}{t^{\alpha\ell_k-\alpha+1+m}}
+ O\left( \frac{1}
{t^{\alpha\ell_k-\alpha+M+2}} \right)\biggr]\biggr\vert
$$
$$
\le \frac{C_6}{(t-t_0)^{\alpha\ell_{K+1}-\alpha+1}}
+ \frac{C(p)}{t^p}  \quad \mbox{for all $p,K,M \in \N$ and $T>T^*$}.  
                                          \eqno{(2.15)}
$$

We can choose $m_1\in \N \cup \{0\}$ such that $\mu_{m_1}:= 
\int^{t_0}_0 (-s)^{m_1}\mu(s) 
ds \ne 0$ and $\mu_j = 0$ for $j=0, 1,2,..., m_1-1$ if $m_1\ge 1$.
Indeed, if such $m_1\in \N\cup \{0\}$ does not exist, then $\mu_m=0$ 
for all $m\in \N \cup \{0\}$.
Therefore Weierstrass polynomial approximation theorem (e.g., Yosida \cite{Y}),
implies $\mu(s)=0$ for almost all $s\in (0,t_0)$, which contradicts the 
assumption $\mu \not\equiv 0$ in $(0,t_0)$.
\\

By (2.15), choosing $M=m_1$, we have
$$
\biggl\vert  
\sum_{k=1}^K \frac{(-1)^{\ell_k+1}}{\Gamma(\alpha-\alpha\ell_k)}
A_k \left(
\begin{array}{c}
-\alpha\ell_k+\alpha-1 \\
m_1 \\
\end{array}
\right)
\frac{\mu_{m_1}}{t^{\alpha\ell_k-\alpha+1+m_1}}
$$
$$
+ O\left( \frac{1}{t^{\alpha\ell_k-\alpha+m_1+2}}\right) \biggr]\biggr\vert
\le \frac{C_6}{(t-t_0)^{\alpha\ell_{K+1}-\alpha+1}}
+ \frac{C(p)}{t^p}, \quad t > T^*.                     \eqno{(2.16)}
$$

In (2.16), we set 
$$
b_k:= 
\left(
\begin{array}{c}
-\alpha\ell_k+\alpha-1 \\
m_1 \\
\end{array}
\right), \quad k\in \N.
$$
We remark that $b_k\ne 0$ because $\alpha\ell_k - \alpha \not\in \N$ 
by the choice of $\ell_k$.

Moreover in (2.16), since $K, p \in \N$ can be arbitrary,
for already fixed $m_1\in \N$,
we can choose sufficiently large $K, p\in \N$ such that
$$
p,\, \alpha_{\ell_{K+1}}-\alpha + 1 > \alpha\ell_1 - \alpha + 1 + m_1.
$$
Then, we can find constants $C_7=C_7(p,K)>0$ and $C_8=C_8(p,K)$ such that 
\begin{align*}
\biggl\vert
&\, \frac{(-1)^{\ell_1+1}}{\Gamma(\alpha-\alpha\ell_1)}A_1
\frac{1}{t^{\alpha\ell_1-\alpha+1+m_1}} 
\left( b_1\mu_{m_1} + \frac{C_7}{t} \right)\\
+&\,  \frac{(-1)^{\ell_2+1}}{\Gamma(\alpha-\alpha\ell_2)}A_2
\frac{1}{t^{\alpha\ell_2-\alpha+1+m_1}} 
\left( b_2\mu_{m_1} + \frac{C_7}{t} \right)\\
+& \, \cdots \cdots \\
+ &\, \frac{(-1)^{\ell_K+1}}{\Gamma(\alpha-\alpha\ell_K)}A_K
\frac{1}{t^{\alpha\ell_K-\alpha+1+m_1}} 
\left( b_K\mu_{m_1} + \frac{C_7}{t} \right)\biggr\vert 
\end{align*}
$$
\le C_8\left( 
\frac{1}{(t-t_0)^{\alpha\ell_{K+1}-\alpha+1}} + \frac{1}{t^p}\right)
\quad \mbox{for each $t>T^*$}.                 \eqno{(2.17)}
$$
The multiplication by $t^{\alpha\ell_1-\alpha+1+m_1}$, yields
\begin{align*}
\biggl\vert
&\frac{(-1)^{\ell_1+1}}{\Gamma(\alpha-\alpha\ell_1)}A_1
\left( b_1\mu_{m_1} + \frac{C_7}{t} \right)
+ \, \frac{(-1)^{\ell_2+1}}{\Gamma(\alpha-\alpha\ell_2)}A_2
\frac{t^{\alpha\ell_1-\alpha+1+m_1}}{t^{\alpha\ell_2-\alpha+1+m_1}} 
\left( b_2\mu_{m_1} + \frac{C_7}{t} \right)\\
+& \, \cdots \cdots 
+ \, \frac{(-1)^{\ell_K+1}}{\Gamma(\alpha-\alpha\ell_K)}A_K
\frac{t^{\alpha\ell_1-\alpha+1+m_1}}{t^{\alpha\ell_K-\alpha+1+m_1}} 
\left( b_K\mu_{m_1} + \frac{C_7}{t} \right)\biggr\vert \\
\le & C_8\left(
\frac{t^{\alpha\ell_1-\alpha+1+m_1}}{(t-t_0)^{\alpha\ell_{K+1}-\alpha+1}}
+ \frac{t^{\alpha\ell_1-\alpha+1+m_1}}{t^p}\right)
\quad \mbox{for $t>T^*$}.
\end{align*}
In view of $\alpha\ell_1 - \alpha + 1 + m_1 
< \alpha\ell_{K+1}-\alpha+1$, $\alpha\ell_1 - \alpha + 1 + m_1 < p$
and $\ell_1 < \ell_2 < \cdots < \ell_K$, 
letting $t \to \infty$, we obtain
$$
\frac{(-1)^{\ell_1+1}}{\Gamma(\alpha-\alpha\ell_1)}A_1b_1\mu_{m_1}
= 0.
$$
By $\mu_{m_1} \ne 0$ and $b_1\ne 0$, we reach $A_1=0$.

Substituting $A_1=0$ into (2.17), we have
\begin{align*}
\biggl\vert
& \frac{(-1)^{\ell_2+1}}{\Gamma(\alpha-\alpha\ell_2)}
\frac{A_2}{t^{\alpha\ell_2-\alpha+1+m_1}} 
\left( b_2\mu_{m_1} + \frac{C_7}{t} \right)\\
+& \, \cdots \cdots \\
+ &\, \frac{(-1)^{\ell_K+1}}{\Gamma(\alpha-\alpha\ell_K)}
\frac{A_K}{t^{\alpha\ell_K-\alpha+1+m_1}} 
\left( b_K\mu_{m_1} + \frac{C_7}{t} \right)\biggr\vert \\
\le &\,  
C_8\left(
\frac{1}{(t-t_0)^{\alpha\ell_{K+1}-\alpha+1}}
+ \frac{1}{t^p}\right)
\quad \mbox{for $t>T^*$}.
\end{align*}
Again we choose larger $p$ and $K\in \N$ such that 
$$
p,\, \alpha\ell_{K+1} - \alpha+1 > \alpha\ell_2 - \alpha + 1 + m_1.
$$
Multipying by $t^{\alpha\ell_2-\alpha+1+m_1}$ and letting $t\to \infty$, we 
obtain
$$
\frac{(-1)^{\ell_2+1}}{\Gamma(\alpha-\alpha\ell_2)}A_2b_2\mu_{m_1}
= 0.
$$

By $\mu_{m_1} \ne 0$ and $b_2\ne 0$, we reach $A_2=0$.

Choosing $p,K \in \N$ at each step and continuing this argument, we finally 
prove that $A_k=0$ for all 
$k\in \N$, that is,
$$
\sumn \frac{a_n}{\la_n^{\ell_k}} = 0, \quad k\in \N.
$$
Thus Lemma 3 completes the proof of Proposition 1.
$\blacksquare$
\section{Proof of Theorem 1}

Assume that for each $p\in \N$, we can find a constant $C(p) > 0$ such that 
(1.7) holds.
We arbitrarily choose $v \in C^{\infty}_0(\Omega)$ satisfying supp $v 
\subset \omega$, and set 
$$
a_n:= (P_nf,\, v)_{L^2(\omega)}, \quad n\in \N.
$$
Then by $\mu=0$ in $(t_0,T_0)$, by (2.1) and the assumption (1.7),
there exists $T_2=T_2(p) > 0$ such 
that 
$$
\left\Vert \sumn \psi_n(t)P_nf\right\Vert_{L^2(\omega)} \le 
\frac{C(p)}{t^p}, \quad t\ge T_2.
$$
Since the series is convergent in $L^2(0,T;L^2(\omega))$, 
taking the scalar products with $v$ in $L^2(\omega)$, we see
$$
\left\vert \sumn a_n\psi_n(t)\right\vert 
\le \frac{C(p)\Vert v\Vert_{L^2(\omega)}}{t^p}, \quad t\ge T_2.
                              \eqno{(3.1)}
$$
\\

Next we prove
$$
\sumn \left\vert \frac{a_n}{\la_n}\right\vert < \infty
                                               \eqno{(3.2)}
$$
for $a_n=(P_nf,\, v)_{L^2(\omega)}$, $n\in \N$.
\\
{\bf Proof of (3.2).}
\\
By $f\in \DDD(A^{\sigma_1})$, we can see $g:= A^{\sigma_1}f
\in L^2(\OOO)$.
Then 
$$
(f, \, \va_{nj}) = (A^{-\sigma_1}g,\, \va_{nj})
= (g,\, A^{-\sigma_1}\va_{nj}) = \la_n^{-\sigma_1}(g,\, \va_{nj}).
$$
Therefore,
$$
P_nf = \sum_{j=1}^{d_n} (f,\, \va_{nj})\va_{nj}
= \la_n^{-\sigma_1}\sum_{j=1}^{d_n} (g,\, \va_{nj})\va_{nj}
= \la_n^{-\sigma_1}P_ng.
$$
Since $v \in C^{\infty}_0(\Omega) \subset \DDD(A)$, we have
$\www{v}:= Av \in L^2(\OOO)$, and so 
$$
a_n = (P_nf, v)_{L^2(\omega)}
= (P_nf, v)_{L^2(\OOO)}
= (\la_n^{-\sigma_1}P_ng, \, A^{-1}\www{v})\\
$$
$$
= (\la_n^{-\sigma_1}A^{-1}P_ng,\, \www{v}) 
= \la_n^{-\sigma_1-2}(P_ng, \, \www{v}).
                                         \eqno{(3.3)}
$$
Therefore,
$$
\sumn \left\vert \frac{a_n}{\la_n}\right\vert 
= \sumn \la_n^{-\sigma_1-2}\vert (P_ng, \www{v})\vert
\le \Vert g\Vert\Vert \www{v}\Vert \sumn \la_n^{-\sigma_1-2}.
$$

On the other hand, by $\mu_n$, $n\in \N$, we renumber the eigenvalues 
$\la_n$ of $A$ repeatedly according to the multiplicities:
$$
\mu_k = \la_1 \,\, \mbox{for $1\le k\le d_1$}, \quad
\mu_k = \la_2 \,\, \mbox{for $d_1+1\le k \le d_1+d_2$}, \cdots.
$$
Then $\mu_n\le \la_n$ for $n\in \N$.

We know that there exists a constant $c_1>0$ such that 
$$
\mu_n = c_1n^{\frac{2}{d}} + o(1) \quad \mbox{as $n \to \infty$}
$$
(e.g., Agmon \cite{Ag}, Theorem 15.1).  Here we recall that 
$d$ is the spatial dimensions.
Therefore, we can find a constant $c_2 > 0$ such that
$$
\la_n \ge c_2n^{\frac{2}{d}} \quad \mbox{as $n \to \infty$}.
                                       \eqno{(3.4)}
$$
Hence, (1.6) yields $\frac{2}{d}(\sigma_1+2) > 1$, so that 
$$
\sumn \la_n^{-\sigma_1-2} \le C_3\sumn n^{-\frac{2}{d}(\sigma_1+2)}
< \infty,
$$
which means (3.2). 
$\blacksquare$

Thus we can apply Proposition 1 in terms of (3.1) and (3.2) to 
obtain $a_n= (P_nf, \, v)_{L^2(\omega)} = 0$ for all $n\in \N$ and 
each $v\in C^{\infty}_0(\omega)$.

Since $v\in C^{\infty}_0(\omega)$ is arbitrary, we obtain
$$
P_nf = 0 \quad \mbox{in $\omega$ for all $n\in \N$}.  \eqno{(3.5)}
$$
For each $n\in \N$, we have $(A-\la_n)P_nf = 0$ in $\OOO$.
In view of (3.5), we can apply the unique continuation for the elliptic
operator $A-\la_n$ (e.g., Isakov \cite{Is}), and we can conclude
that $P_nf = 0$ in $\OOO$ for each $n\in \N$.
Hence $f = \sumn P_nf = 0$ in $L^2(\OOO)$.
Thus the proof of Theorem 1 (i) is complete.

Finally in the case of data $\NUNU  u\vert_{\gamma\times (t_0,T_)}$,
as test functions we choose 
$v \in C^{\infty}(\gamma)$ and set 
$a_n= (\ppp_{\nu_A}P_nf,\, v)_{L^2(\gamma)}$ for $n\in \N$.
Then for any $\ep>0$, we have
$$
\Vert \NUNU P_nf\Vert_{L^2(\gamma)} 
\le C\Vert P_nf\Vert_{H^{\frac{3}{2}+2\ep}(\OOO)}
\le C\Vert A^{\frac{3}{4}+\ep}P_nf\Vert 
\le C\la_n^{\frac{3}{4}+\ep}\Vert P_nf\Vert
$$
by the trace theorem.  

By $f \in \DDD(A^{\sigma_2})$ where $\sigma_2>0$ is given by (1.6), we have
$A^{\sigma_2}f=: g \in L^2(\OOO)$, and so 
$$
(f, \, \va_{nj}) = (A^{-\sigma_2}g,\, \va_{nj})
= (g,\, A^{-\sigma_2}\va_{nj}) = \la_n^{-\sigma_2}(g,\, \va_{nj}).
$$
Therefore, we can verify
$$
P_nf = \sum_{j=1}^{d_n} (f,\, \va_{nj})\va_{nj}
= \la_n^{-\sigma_2}\sum_{j=1}^{d_n} (g,\, \va_{nj})\va_{nj}
= \la_n^{-\sigma_2}P_ng.                                
$$
Therefore, 
$$
\vert a_n\vert \le \Vert \ppp_{\nu_A}P_nf\Vert_{L^2(\gamma)}
\Vert v\Vert_{L^2(\gamma)}
\le C_4\la_n^{\frac{3}{4}+\ep}\Vert P_nf\Vert 
\le C_4\la_n^{\frac{3}{4}+\ep-\sigma_2}\Vert P_ng\Vert,
$$
and so (3.4) yields
$$
\sumn \left\vert \frac{a_n}{\la_n}\right\vert 
\le C_4\sumn \la_n^{-\frac{1}{4}+\ep-\sigma_2}\Vert P_ng\Vert
\le C_4\sumn n^{-(\frac{1}{4}-\ep+\sigma_2)\frac{2}{d}}\Vert g\Vert.
                                      \eqno{(3.6)}
$$
In terms of (1.6), for sufficiently small $\ep>0$ we see that 
$\left( \frac{1}{4}-\ep +\sigma_2\right)\frac{2}{d} > 1$ and 
(3.2) is verified to hold in the case of 
$a_n = (\ppp_{\nu_A}P_nf,\, v)_{L^2(\gamma)}$.

Consequently we can 
argue similarly to the case of $u\vert_{\omega\times (t_0,T_0)}$
to complete the proof of Theorem 1 (ii).
$\blacksquare$
\section{Proof of Corollary 1.}

{\bf First Step: holomorphicity of time factors.}
\\
We consider the holomorphicity of related functions to (2.1).

For $\theta \in \left(0, \frac{\pi}{2}\right)$ and $z_0\in \C$, we define
a sector $S_{\theta}(z_0) \subset \C$ by 
$$
S_{\theta}(z_0):= \{ z\in \C;\, z\ne z_0. \quad
\vert \mbox{arg}\, (z-z_0) \vert < \theta\}.
$$
Let $t_0 < T_1 < T_0$.
Then
\\
{\bf Lemma 4.}
\\
{\it
Let $h = h(z)$ be holomorphic in $S_{\theta}(T_1-t_0)$ and let $\mu \in L^1(0,t_0)$.Then
$$
H(z):= \int^{t_0}_0 h(z-s)\mu(s) ds
$$
is holomorphic in $z\in S_{\theta}(T_1)$.
}
\\
{\bf Proof of Lemma 4.}
\\
For $z \in S_{\theta}(T_1)$, we choose $r>0$ and $0<\theta_1<\theta$ such that 
$z = T_1 + re^{\sqrt{-1}\theta_1}$.  Therefore,
$z-s = (T_1-s) + re^{\sqrt{-1}\theta_1} \in S_{\theta_1}(T_1-s)
\subset S_{\theta_1}(T_1-t_0)$ by $T_1-t_0 \le T_1-s$  
for all $s \in [0,t_0]$.  Hence, $z-s \in S_{\theta_1}(T_1-t_0)$ if
$z\in S_{\theta}(T_1)$ and $s\in [0,t_0]$.

Consequently, by the Lebesgue convergence theorem and $\mu \in L^1(0,T)$, 
we have
$$
\frac{\ppp H}{\ppp z}(z) = \int^{t_0}_0 \frac{\ppp h}{\ppp z}(z-s)
\mu(s) ds, \quad z\in S_{\theta}(T_1).
$$
Thus the holomorphicity follows.
$\blacksquare$
\\

We see that $g(z):= z^{\alpha-1}E_{\alpha,\alpha}(-\la_nz^{\alpha})$ is 
holomorphic in $z \in S_{\theta}(T_1-t_0)$ by $T_1-t_0 > 0$.
We can choose a small constant $\theta_0>0$ such that 
$$
\vert E_{\alpha,\alpha}(-\la_n(z-s)^{\alpha}) \vert 
\le \frac{C}{1+\la_n\vert z-s\vert^{\alpha}}         \eqno{(4.1)}
$$
for $\vert \mbox{arg}\, (z-s)\vert \le \theta_0$. 
\\
{\bf Proof of (4.1).}
\\
We choose constants $\mu_1, \mu_2$ satisfying
$$
\left\{ \begin{array}{rl}
&\frac{\pi}{2}\alpha < \mu_1 < \pi\alpha \quad \mbox{if $0<\alpha<1$},\\
&\frac{\pi}{2}\alpha < \mu_2 < \pi \quad \mbox{if $1<\alpha<2$}.
\end{array}\right.
$$
Setting $\eta := z-s$, we apply Theorem 1.6 (p.35) in \cite{Po}, so 
that we know that (4.1) holds if 
$$
\mu_1 \le \vert \mbox{arg}\,(-\la_n\eta^{\alpha})\vert \le \pi
\quad \mbox{if $0<\alpha<1$}
$$
and
$$
\mu_2 \le \vert \mbox{arg}\,(-\la_n\eta^{\alpha})\vert \le \pi
\quad \mbox{if $1<\alpha<2$}.
$$
Henceforth, for $0\le \theta < \pi$, we use that  
$\theta \le \vert \mbox{arg}\, (-\xi)\vert \le \pi$ is equivalent to  
$\vert \mbox{arg}\,\xi \vert \le \pi - \theta$.

Let $0<\alpha<1$.  Then, $\vert \mbox{arg}\,\eta^{\alpha}\vert \le 
\pi - \mu_1$ implies $\mu_1 \le \vert \mbox{arg}\,(-\eta^{\alpha})\vert
\le \pi$.  Setting $\theta_0:= \pi - \mu_1$, we see that if 
$\vert \mbox{arg}\, \eta\vert \le \theta_0$, then 
$$
\vert \mbox{arg}\,\eta^{\alpha}\vert = \alpha \vert \mbox{arg}\,\eta\vert
\le \alpha\theta_0 = \alpha(\pi - \mu_1) < \pi - \mu_1,
$$
so that (4.1) holds for $\vert \mbox{arg}\,(z-s)\vert \le \theta_0$
in the case of $0<\alpha<1$.

Next let $1<\alpha<2$.  We choose $\mu_2 = \pi - \ep_0$ with sufficiently 
small constant $\ep_0>0$.  
Let $\vert \mbox{arg}\,\eta \vert \le \ep_0$.  Then
$\vert \mbox{arg}\,\eta^{\alpha}\vert \le \alpha\ep_0$, and so
$$
\pi - \alpha\ep_0 \le \vert \mbox{arg}\,(-\eta^{\alpha})\vert \le \pi.
$$
For $\mu_2 \in \left(\frac{\pi}{2}\alpha,\, \pi\right)$, we further choose
small $\ep_0>0$ such that $\mu_2 < \pi - \alpha\ep_0$.
Then we obtain
$\mu_2 < \vert \mbox{arg}\,(-\eta^{\alpha})\vert \le \pi$, that is,
(4.1) holds for $\vert \mbox{arg}\,(z-s)\vert \le \ep_0$
in the case of $1<\alpha<2$.

Now, by choosing 
$$
\theta_0 = 
\left\{ \begin{array}{rl}
\pi - \mu_1, \quad & 0<\alpha<1,\\
\ep_0, \quad & 1< \alpha < 2,
\end{array}\right.
$$
the estimate (4.1) holds.
$\blacksquare$

Therefore, since $z \in S_{\theta_0}(T_1)$ implies that 
$\vert \mbox{arg}\,(z-s) \vert \le \theta_0$ for any 
$s \in [0,t_0]$ by $s \le t_0 < T_1$, we have
$$
\vert E_{\alpha,\alpha}(-\la_n(z-s)^{\alpha})\vert 
\le \frac{C}{1 + \la_n\vert z-s\vert^{\alpha}}\quad
\mbox{for all $z \in S_{\theta_0}(T_1)$, $s\in [0,t_0]$ 
and all $n\in \N$.}                                       \eqno{(4.2)}
$$
Moreover, for arbitrarily fixed $s \in [0,t_0]$, the function in 
$z$
$$
g(z):= z^{\alpha-1} E_{\alpha,\alpha}(-\la_nz^{\alpha})
$$
is holomorphic in $S_{\theta_0}(T_1-t_0)$, because $T_1-t_0 > 0$.

{\bf Second Step.}
\\
We can assume that $\mu \not\equiv 0$ in $(0,t_0)$.
Indeed if $\mu\equiv 0$ in $(0,t_0)$, the the proof is already 
complete.

From Lemma 4 it follows that $\psi_n(t)$, $n\in \N$
defined by (2.2) for 
$t\ge T_1$, can be extended to a holomorphic function 
$\psi_n(z)$ for $z \in S_{\theta_0}(T_1)$.

Assume that $ u(x,t) = 0$ for $x\in \omega$ and $t_0<t < T_0$.
We will prove
$$
\sumn \sup_{z\in S_{\theta_0}(T_1)} \vert \psi_n(z)\vert \vert a_n\vert
< \infty.                                      \eqno{(4.3)}
$$
Here we set $a_n=(P_nf, v)_{L^2(\omega)}$, $n\in \N$ for 
$v\in C^{\infty}_0(\omega)$ and 
$a_n = (\ppp_{\nu_A}P_nf,\, v)_{L^2(\gamma)}$, $n\in \N$ for 
$v\in C^{\infty}(\gamma)$.
\\
{\bf Proof of (4.3).}
\\
Since $\vert z-s\vert \ge T_1-t_0> 0$ for each $z \in S_{\theta_0}(T_1)$ and 
$0\le s \le t_0$, we estimate $\vert \psi_n(z)\vert$. 
In view of (4.2), similarly to Lemma 2, we have
\begin{align*}
& \vert \psi_n(z)\vert 
= \left\vert \int^{t_0}_0 (z-s)^{\alpha-1}
E_{\alpha,\alpha}(-\la_n(z-s)^{\alpha}) \mu(s) ds\right\vert \\
\le &C\int^{t_0}_0 \frac{\vert z-s\vert^{\alpha-1}}
{1 + \la_n\vert z-s\vert^{\alpha}} \vert \mu(s)\vert ds
\le \frac{C}{\la_n} \Vert \mu\Vert_{L^1(0,t_0)},
\quad z\in S_{\theta_0}(T_1), \, n\in \N. 
\end{align*}
Hence, applying (3.2) and (3.6) respectively for $a_n=(P_nf,v)_{L^2(\omega)}$
and $a_n=(\ppp_{\nu_A}P_nf, \, v)_{L^2(\gamma)}$, we obtain
$$
\sumn \sup_{z\in S_{\theta_0}(T_1)} \vert \psi_n(z)\vert
\vert a_n\vert
\le \sumn \frac{C}{\la_n}\Vert \mu\Vert_{L^1(0,t_0)}
\vert a_n\vert < \infty.
$$
Thus the proof of (4.3) is complete.
$\blacksquare$

By means of (4.3) and the holomorphicity of $\psi_n(z)$, we see that 
$\sumn a_n\psi_n(z)$ is holomorphic in any compact set in 
$S_{\theta_0}(T_1)$, and in particular, $\sumn a_n\psi_n(t)$ is real
analytic in $t>T_1$.  Consequently if $u=0$ in $\omega \times 
(t_0,T_0)$ or $\ppp_{\nu_A}u = 0$ on $\gamma \times (t_0,T_0)$, then
$$
\sumn a_n\psi_n(t) = 0, \quad t>T_1.
$$
Since $v$ is arbitrarily chosen, we see that $u=0$ in $\omega \times 
(t_0,\infty)$ or $\ppp_{\nu_A}u = 0$ on $\gamma \times (t_0,\infty)$.
Thus the assumptions in Theorem 1 are satisfied, so that 
$f=0$ in $\OOO$ follows, and the proof of Corollary 1 is complete.
$\blacksquare$
\section{Concluding remarks}

{\bf 5.1.}
The proof of the main result relies on Proposition 1 which is concerned with 
the infinite series of $\psi_n(t)$, $n\in \N$ defined by (2.2).
The proof of the proposition is based on the asymptotic expansion of 
$\int^{t_0}_0 \frac{\mu(s)}{(t-s)^{\sigma}} ds$ with 
$\sigma = \alpha k - \alpha+1$.
Here in the case of $0<\alpha<1$, we extract its essence and prove as 
\\
{\bf Proposition 2.}
\\
{\it
(i) We assume that there exist $a\in \R$ and $t_0 > 0$ and 
$v-a \in H_{\alpha}(0,T)$ satisfies 
$$
\pppa (v-a)(t) = 0, \quad t>t_0                 \eqno{(5.1)} 
$$
and for each $p\in \N$ we can find a constant $C(p)>0$ such that 
$$
\vert v(t)\vert \le \frac{C(p)}{t^p} \quad \mbox{as $t \to \infty$}.
                                                \eqno{(5.2)}
$$
Then $a=0$ and $v(t)=0$ for $t>0$. 
\\
(ii) In particular, if $\pppa (v-a)(t) = v(t) = 0$ for $t>t_0$
with some $a\in \R$ and $t_0>0$, then $v(t) = 0$ for $t>0$.
}
\\

Proposition 2 (ii) corresponds to Theorem 1 in Kinash and Janno \cite{KJ}, 
which treats more general time-fractional derivative.  We can relax (5.1) as 
the vanishing in some finite interval $t_0<t<T_0$ but we do not discuss.
\\
{\bf Proof of Proposition 2 (i).}
\\
Part (ii) is a direct consequence of part (i), and so it suffices to prove
part (i).
Setting $\mu(t) := \pppa (v-a)$ for $t>0$, we have $\mu \in L^2(0,T)$ by 
$v-a \in H_{\alpha}(0,T)$ and 
$$
v(t) = a + \frac{1}{\Gamma(\alpha)}\int^t_0 (t-s)^{\alpha-1}\mu(s) ds,
\quad t>0
$$
(e.g., \cite{KRY}, \cite{Y2022}).
By (5.1), we obtain
$$
v(t) = a + \frac{1}{\Gamma(\alpha)}\int^{t_0}_0 (t-s)^{\alpha-1}\mu(s) ds,
\quad t>t_0.                                      \eqno{(5.3)}
$$

Let $T_1 > t_0$ be arbitrarily chosen.  Applying (2.14), for any $M\in \N$ we
have
$$
\int^{t_0}_0 \frac{\mu(s)}{(t-s)^{1-\alpha}} ds
= \sum_{m=0}^M 
\left(
\begin{array}{c}
1-\alpha \\
m \\
\end{array}
\right) \frac{\mu_m}{t^{1-\alpha+m}}
+ O\left( \frac{1}{t^{1-\alpha+M+1}}\right), \quad t> T_1.
$$
Here we set $\mu_m := \int^{t_0}_0 (-s)^m\mu(s) ds$ for 
$m\in \N \cup \{ 0\}$.
Substituting this into (5.3), we reach 
$$
v(t) = a + \sum_{m=0}^M \frac{1}{\Gamma(\alpha)} 
\left(
\begin{array}{c}
1-\alpha \\
m \\
\end{array}
\right) \frac{\mu_m}{t^{1-\alpha+m}}
+ O\left( \frac{1}{t^{1-\alpha+M+1}}\right), \quad t> T_1.
$$
We number the set $\{ m\in \N;\, \mu_m \ne 0\}$ as 
$$
m_1 < m_2 < \cdots.
$$
Henceforth we assume that such a set of $\mu_k$ is an infinite set.
The proof is same for the case of the finite set.

Applying (5.2) for all $p\in \N$, we obtain
$$
\left\vert a + \sum_{j=1}^{m_N} \frac{1}{\Gamma(\alpha)} 
\left(
\begin{array}{c}
1-\alpha \\
m \\
\end{array}
\right) \frac{\mu_{m_j}}{t^{1-\alpha+m_j}}
+ O\left( \frac{1}{t^{1-\alpha+m_N+1}} \right)\right\vert 
\le \frac{C(p)}{t^p}
                                          \eqno{(5.4)}
$$
as $t\to \infty$. 
Letting $t\to \infty$ in (5.4), we obtain $a=0$.

For arbitrarily fixed $N\in \N$, we can choose
$p > 1-\alpha+m_N+1$.  
Multiplying (5.4) with $t^{1-\alpha+m_1}$ and letting $t\to \infty$, we obtain
$$
\frac{1}{\Gamma(\alpha)} 
\left(
\begin{array}{c}
1-\alpha \\
m_1 \\
\end{array}
\right) \mu_{m_1} = 0,
$$
that is, $\mu_{m_1} = 0$.  Continuing this argument, we reach
$\mu_{m_1} = \mu_{m_2} = \cdots = \mu_{m_N} = 0$.
Since $N\in \N$ is arbitrary, we see that $\mu_{m_k} = 0$ for $k\in \N$.
The Weierstrass polynomial approximation theorem implies $\mu = 0$ in 
$(0,t_0)$.  Thus $a=0$ and $v=0$ in $(0,t_0)$ by (5.3).
$\blacksquare$

{\bf 5.2.}
We sum up the uniqueness results for the inverse source problem of determining
a spatially varying function $f(x)$ in (1.5):
\begin{itemize}
\item
{\bf Data over time interval $0<t<T_0$:}\\
This is a conventional formulation and there are satisfactory researches 
on the uniqueness and also the conditional stability.
\item
{\bf Data over time interval $0<t_0 < t <T_0$:}\\
The current work proved the uniqueness if $\mu \not\equiv 0$ in $(0,t_0)$ and 
$\mu=0$ in $(t_0,\,T_0)$.
\\
The uniqueness remains open under other condition  
$\mu \not\equiv 0$ in $(t_0,\, T_0)$, in particular, 
$\{ t\in (t_0,\, T_0);\, \mu(t) = 0\}$ has no interior points.
\end{itemize}

{\bf 5.3.}
In this article, for keeping conciseness, we are restricted to establish the 
uniqueness of $f(x)$ or $\mu(t)$ by fixing the order $\alpha \in (0,1) \cup 
(1,2)$.  
Indeed our data concerning $u\vert_{\omega\times (t_0,T_0)}$ or
$\NUNU u\vert_{\gamma\times (t_0,T_0)}$ can contain enough 
information for determining also an order $\alpha\in (0,2)$ 
including $\alpha=1$. The method in Yamamoto \cite{Y5}, \cite{Y6} is 
applicable, but we postpone the details to a forthcoming work.

The order $\alpha \in (0,2)$ is an important parameter characterizing 
the anomaly of the diffusion, which is a quantity describing the heterogeneity 
of media under consideration.  Usually also $\alpha$ cannot be 
determined a priori and we should estimate or identify $\alpha$ by available 
observation data.  Thus the inverse problems of determining an order $\alpha$ 
are significant.
 
{\bf 5.4.}
In this article, we assume that initial values are zero.
On the other hand, our data under assumption 
$\mu(t) = 0$ for $t>t_0$ may determine also an initial
value uniquely and in a forthcoming paper we will discuss.

{\bf 5.5.}
Our main result asserts that if  
$$
\mu \not\equiv 0 \quad \mbox{in $(0,\, t_0)$},               \eqno{(5.5)}
$$
then the uniqueness in determining $f$ in $\OOO$ follows for 
$\alpha \in (0,1) \cup (1,2)$.
As is seen from 
Theorem 4 in \cite{CLY}, for the case of $\alpha=1$, 
the condition (5.5) is not at all sufficient for the uniqueness 
for the corresponding 
inverse source problem, and we can prove the uniqueness under 
much stronger condition 
$$
\int^{t_0}_0 e^{\la_ns}\mu(s) ds \ne 0 \quad \mbox{for all $n\in \N$}.
                                    \eqno{(5.6)}
$$
This difference in the conditions for the uniqueness,
can be interpreted by a kind of anomaly of the fractional
diffusion as that the case $\alpha \ne 1$ keeps more information of the state 
before the observation time $(t_0, \, T_0)$, while in the case of $\alpha=1$, 
we need more data holding information of every eigenmode 
corresponding to the eigenvalue $\la_n$, $n\in \N$.
\\

{\bf Acknowledgments.}
The author was supported by Grant-in-Aid for Scientific Research
Grant-in-Aid (A) 20H00117 of 
Japan Society for the Promotion of Science.


\end{document}